\documentclass{article}
\usepackage{amssymb,amsmath,amscd,amsthm}
\begin{document}

\newtheorem{theorem}{Theorem}
\newtheorem{lemma}{Lemma}

{\centerline {\Large {\bf  Elementary 
equivalence of the semigroup of invertible}}}

{\centerline {\Large {\bf  matrices with nonnegative elements}}}
\vspace{.3truecm}

\rightline{\large E.I. Bunina, A.V. Mikhalev}

\vspace{.2truecm}

Let $R$ be a linearly ordered ring with $1/2$, $G_n(R)$ ($n\ge 3$) be the
subsemigroup of $GL_n(R)$ consisting of all
matrices with nonnegative elements. In~[1],
there is a description of all automorphisms
of the semigroup $G_n(R)$ in the case where  $R$ is a skewfield and $n\ge 2$. 
In~[2], there is a description of all automorphisms of
the semigroup $G_n(R)$ for the case where
$R$ is an arbitrary linearly ordered ring with $1/2$ and $n\ge 3$.
In this paper, we classify semigroups $G_n(R)$ up to elementary equivalence.

Two models ${\cal U}_1$ and ${\cal U}_2$ of the same first order 
language~$\cal L$
(e.\,g. two groups, semigroups, or two rings, semirings) are called
\emph{elementarily equivalent}, if every sentence~$\varphi$ of the language~$\cal L$
is true in the model~${\cal U}_1$ if and only if
it is true in the model~${\cal U}_2$. 

Any two finite models of the same language
are elementarily equivalent if and only if they are isomorphic.
Any two isomorphic models are elementarily equivalent but for
infinite models the converse is not true. For example, the field $\mathbb C$
of complex numbers and the field $\overline {\mathbb Q}$ of algebraic numbers
are elementarily equivalent but not isomorphic (since they have different cardinalities).

The first results in elementary equivalence of linear groups were obtained
by A.I. Maltsev in~[3].
He proved that the groups ${\cal G}_m(K)$ and ${\cal G}_n(K')$
(where ${\cal G}=GL,PGL,SL,PSL$, $m,n> 2$, $K$ and $K'$ are fields  of characteristic~$0$)
are elementarily equivalent if and only if $m=n$
and the fields $K$ and $K'$ are elementarily equivalent.

In 1992, C.I.~Beidar and A.V.~Mikhalev ([4]), 
using some results of model theory
 (namely, the construction of ultrapower and  Keisler--Shelah Isomorphism Theorem)
formulated a general approach to the problem of elementary equivalence
of various algebraic structures.

Taking into account some results of the theory of linear groups over rings,
they obtained easy proofs of theorems similar to Maltsev's theorem
in rather general situations (for linear groups over prime rings, 
multiplicative semigroups, lattices of submodules, and so on).

In 1998--2004, E.I.~Bunina continued to study elementary properties
of linear groups. In 1998, the results of A.I.~Maltsev were
generalized to unitary linear groups over fields with involutions ([5]), 
and then to unitary linear groups over rings an skewfields with involutions ([6]).
In 2001--2004  ([7])  similar results were obtained for Chevalley groups
over fields.

We use  notations and definition from~[2].

Now we will recall the most necessary definitions.

Suppose that $R$ is a linearly ordered ring, $R_+$ is
the set of all positive elements, $R_+\cup \{ 0\}$ is the set
of all nonnegative elements of the ring~$R$. By $G_n(R)$
we denote the subsemigroup of $GL_n(R)$ consisting of all
matrices with nonnegative elements.

The set of all invertible elements of~$R$ is denoted by~$R^*$.
The set $R_+\cap R^*$ is denoted by~$R_+^*$.
If $T\subset R$, then $Z(T)$ denotes the center of~$T$,
$Z^*(T)=Z(T)\cap R^*$, $Z_+(T)=Z(T)\cap R_+$, $Z_+^*(T)=Z(T)\cap R_+^*$.

Let 

$I=I_n$,

$\Gamma_n(R)$ be the group consisting of all invertible matrices from $G_n(R)$,

$\Sigma_n$ be the symmetric group of  order~$n$,

$S_\sigma$ be the matrix of the permutation
 $\sigma\in \Sigma_n$ (i.\,e., the matrix
$(\delta_{i\sigma(j)})$),

$S_n=\{ S_\sigma|\sigma\in \Sigma_n\}$,

$diag[d_1,\dots,d_n]$ be the diagonal matrix with elements
 $d_1,\dots,d_n$ on the diagonal, where
$d_1,\dots,d_n\in R_+^*$,

$D_n(R)$ be the group of all invertible diagonal matrices from $G_n(R)$,

$D_n^Z(R)$ be the center of $D_n(R)$.

By $K_n(R)$ we denote the subsemigroup in $G_n(R)$ consisting of all matrices
of the form
$$
\begin{pmatrix}
X_{n-1}& 0\\
0& x
\end{pmatrix},\quad X_{n-1}\in G_{n-1}(R),\ x\in R_+^*.
$$

Let $E_{ij}$ be the matrix with one nonzero element $e_{ij}=1$.

By $B_{ij}(x)$ we denote the matrix $I+xE_{ij}$, where $x\in R_+$, $i\ne j$.

In this paper we prove the following theorem:

\begin{theorem}
Semigroups $G_n(R)$ and  $G_m(S)$ $($where $n,m\ge 3$, $1/2\in R$, $1/2\in S)$
are elementarily equivalent if and only if $n=m$ and the semirings
$R_+$ and $S_+$ are elementarily equivalent.
\end{theorem}

\begin{lemma}\label{l1}
The formula
$$
Invert(M):= \exists X(MX=XM=1)
$$
is true in the semigroup $G_n(R)$ for elements of the group    $\Gamma_n(R)$,
and only for them.
\end{lemma}
The proof is obvious.

\begin{lemma}\label{l2}
The formula 
$$
Inv(M):=\bigl( (M^2=1)\land M\ne 1\bigr)
$$
is true in the semigroup $G_n(R)$ for involutions, and only for them.
These involutions have the form $diag[t_1,\dots,t_n]S_\sigma$, 
where $\sigma^2=1$,
$t_i\cdot t_{\sigma(i)}=1$ for all $i=1,\dots,n$.
\end{lemma}
\begin{proof}
The first part is clear; the second part follows from Lemma~3 ([2]).
\end{proof}

\begin{lemma}\label{l3}
\emph{(1)} There exists a formula $Diag(M)$ that is true in the
semigroup $G_n(R)$ for the matrices $M\in D_n(R)$, and only for them.

\emph{2)} There exists a formula $CDiag(M)$ that is true in the
semigroup $G_n(R)$ for the matrices $M\in D_n^Z(R)$, and only for them.
\end{lemma}
\begin{proof}
Consider the formula
$$
ComCon(A):= \bigl( Invert(A)\land \forall M(\exists N (M=NAN^{-1})\Rightarrow
AM=MA)\bigr).
$$
This formula states that a matrix~$A$ is invertible and it commutes
with all matrices conjugate to~$A$.
By Lemma~4 ([2]) the formula $ComCon(A)$ is true for all
 $A\in D_n^Z(R)$, it can be true for some elements 
from $D_n(R)\setminus D_n^Z(R)$, and it can not be true for elements 
from $\Gamma_n(R)\setminus D_n(R)$. 

Now let us introduce the new formula:
$$
NComInv(A):= \bigl( ComCon(A)\land \forall M (Inv(M)\Rightarrow AM\ne MA)\bigr).
$$
This formula gives us an additional condition ``a matrix $A$ 
does not commute with any involution''. It follows from the proof of
Lemma~4 ([2]) that this formula is always true for matrices $M\in D_n^Z(R)$
having different eigenvalues, it can be true for some matrices $M\in D_n(R)
\setminus D_n^Z(R)$ having different eigenvalues, and it cannot be true
for any other matrices.

Consider now the formula
$$
Diag(M):= \exists A(NComInv(A)\land MA=AM\land Invert(M)).
$$
The formula $Diag(M)$ is true for matrices from ${\cal X}=D_n(R)$ (see 
the proof of Lemma~4 ([2])), i.\,e., it satisfies  assertion (1) of
the lemma.

The formula 
$$
CDiag(M): = \bigl( Invert(M)\land \forall A(Diag(A)\Rightarrow AM=MA)\bigr)
$$
satisfies the assertion (2) of the lemma.
\end{proof}

\begin{lemma}\label{l4}
For any $n\ge 2$ there exists a sentence $Size_n$ that is true
in all semigroups $G_n(R)$, where $1/2\in R$, and is false in all semigroups
$G_m(S)$, where $m\ne n$, $1/2\in S$.
\end{lemma}
\begin{proof}
Consider the sentence
\begin{multline*}
Size_n:= \exists X_1\dots \exists X_{n!} (\forall M (Diag(M)\Rightarrow
\bigwedge_{i,j=1; i\ne j}^{n!} X_i\ne MX_j)\land\\
\land \forall X(Invert (X)\Rightarrow \exists M(Diag (M)\land \bigvee_{i=1}^{n!}
(X=MX_i)))).
\end{multline*}
It is clear that is satisfies our conditions.
\end{proof}

Now we can suppose that the dimension~$n$ of the semigroup $G_n(R)$
is known, i.e., it  can be used in formulas.

\begin{lemma}\label{l5}
There exists a formula $CDOneMany_n(M)$ that is true for all matrices
$$
diag[\alpha,\dots,\alpha,\beta,\alpha,\dots,\alpha]\in D_n^Z(R),\quad \alpha
\ne \beta,
$$
and only for them.
\end{lemma}
\begin{proof}
From Lemma~5 ([2]) it follows that our matrices~$M$ 
can be characterized by the condition
$$
C_{\Gamma_n(R)}(M)/D_n(R)\cong \Sigma_{n-1}.
$$
Let elements $\sigma_1,\dots, \sigma_N$ (where $N=(n-1)!$) of $\Sigma_{n-1}$
be numerated and $\sigma_i\cdot \sigma_j=\sigma_{\gamma(i,j)}$. Then
we obtain the formula
\begin{multline*}
CDOneMany_n(M): = \exists X_1\dots \exists X_N \Biggl(\left(\bigwedge_{i=1}^N 
X_iM=MX_i\right)\land\\
\land \left(\bigwedge_{i\ne j, i,j=1}^N \forall A(Diag(A)\Rightarrow X_i\ne X_j A
)\right)\land
\left(\bigwedge_{i=1}^N Invert(X_i)\right)\land\\
\land \forall X\left(Invert (X)\land XM=MX\Rightarrow \exists A \left(Diag(A)\land 
\bigvee_{i=1}^N (X=AX_i)\right)\right)\land\\
\land \left(\bigwedge_{i,j=1}^N \exists A
(Diag(A)\land X_i\cdot X_j=AX_{\gamma(i,j)})\right)\Biggr).
\end{multline*}
The lemma is proved.
\end{proof}

Similarly, we can introduce the formula $CDAll_n(M)$, characterizing matrices
$\alpha I\in D_n^Z(R)$.

\begin{lemma}\label{l6}
For any matrix 
$$
M=diag[\alpha,\dots,\alpha,\underbrace{\beta}_i,
\alpha,\dots,\alpha]\in D_n^Z(R),
\quad \alpha\ne \beta,
$$
the formula
\begin{multline*}
DSame_{1,n-1}(A,M):= CDOneMany_n(A)\land 
CDOneMany_n
(M)\land\\
\land (CDOneMany_n(AM)\lor CDAll_n(AM))\land 
(CDOneMany_n(AM^{-1})\lor CDAll_n(AM^{-1}))
\end{multline*}
characterizes matrices
$$
A=diag[\gamma,\dots,\gamma,\underbrace{\delta}_i,\gamma,\dots,\gamma]
\in D_n^Z(R),
\quad \gamma\ne \delta.
$$
\end{lemma}

\begin{proof}
If
$$
A=diag[\gamma,\dots,\gamma,\underbrace{\delta}_i,\gamma,\dots,\gamma]
\in D_n^Z(R),
\quad \gamma\ne \delta,
$$
then by Lemma~\ref{l5} the following holds: 
the subformula  $DOneMany_n(A)$ is true;
$$
AM=diag [\alpha \gamma,\dots,\alpha\gamma,\beta\delta,\alpha\gamma,\dots,
\alpha\gamma]\in D_n^Z(R);
$$
 if $\alpha\gamma=\beta\delta$, then
the subformula $CDAll_n(AM)$ is true; and if $\alpha\gamma\ne
\beta\delta$, then the subformula $CDOneMany_n(AM)$ is true.
Similarly for $AM^{-1}$.

Conversely, let for 
$$
M=diag [\alpha,\dots,\alpha,\beta,\alpha,\dots,\alpha]
$$
and some~$A$ the formula $DSame_{1,n-1}(A,M)$ be true.
Since in this case the subformula $CDOneMany_n(A)$ is true, 
we have 
$$A=diag[\gamma,\dots,\gamma,\delta,\gamma,\dots,\gamma]\in D_n^Z(R), 
quad \delta\ne \gamma.
$$
If $\beta$ and $\delta$ are staying at the same $i$-th place, 
then everything is proved.
Suppose that it is not true. Then, without loss of
generality, we can assume that
$M=diag[\beta,\alpha,\dots,\alpha]$, 
 $A=diag[\gamma,\delta,\gamma,\dots,\gamma]$.
In this case
 $MA=diag[\beta\gamma,\alpha\delta,\gamma\alpha,\dots,\gamma\alpha]$.
It is clear that $\beta\gamma\ne \alpha\gamma$ and
 $\alpha\delta \ne \alpha\gamma$. 
If  $n> 3$, then we immediately get $\neg (CDOneMany_n(AM)\lor 
CDAll_n(AM))$, which contradicts to our assumption.
If $n=3$, then it may be that case that
 $\beta\gamma=\alpha\delta$, and then the formula
$CDOneMany_n(AM)$ is true. In this case, let us consider 
the matrix 
$AM^{-1}=diag[\beta\gamma^{-1},\alpha\delta^{-1},\alpha\gamma^{-1}]$.
Since $\beta\gamma^{-1}\ne \alpha \gamma^{-1}$ and $\alpha \delta^{-1}\ne
\alpha \gamma^{-1}$, we have $\beta\gamma^{-1}=\alpha\delta^{-1}$.
Consequently,
$$
\begin{cases}
\beta\gamma=\alpha\delta\\
\beta\gamma^{-1}=\alpha\delta^{-1}
\end{cases}
\Rightarrow \alpha=\beta,
$$
and this contradicts  our assumption.
\end{proof}

\begin{lemma}\label{l7}
There exists a formula $KOneMany_n(X,M)$ with to free variables
such that for every matrix~$A$ satisfying
the formula $CDOneMany_n(A)$, the set of all matrices~$M$ satisfying
the formula $KOneMany_n(A,M)$ is a group
$\Phi_N (K_n(R))=N(K_n(R))N^{-1}$ for some matrix $N\in \Gamma_n(R)$.
\end{lemma}
\begin{proof}
Consider the  formula
$$
KOneMany_n(X,M):=\exists M' (XM'=M'X\land DSame_{1,n-1}(M',M)).
$$
This formula states that
$$
X\in \bigcup_{A=diag[\alpha,\dots,\alpha,\beta,\alpha,\dots,\alpha]\in D_n^Z(R),
\alpha\ne \beta} C_{G_n(R)}(A),
$$
whence by Lemma~5 ([2]) we obtain our statement.
\end{proof}

\begin{lemma}\label{l8}
There exist formulas $Cycle_n(X)$, $Trans_n(X,Y)$, and $Perm_n(X,Y,Z)$ such
that for any matrices $M_1,M_2,M_3$ satisfying the formulas
$Cycle_n(M_1)$, $Trans_n(M_1,M_2)$ and $Perm_n(M_1,M_2,M_3)$, 
there exists a matrix $N\in \Gamma_n(R)$ such that
$M_1=\Phi_N(S_{(1,2,\dots,n)})$, $M_2=\Phi_N(S_{(1,2)})$, 
$M_3=\Phi_N(S_\sigma)$
for some $\sigma\in \Sigma_n$.
\end{lemma}

\begin{proof}
Consider some matrix $M$ satisfying the formula
$$
Cycle_n(M):= (M^n=1)\land \forall X (CDiag(X)\land MX=XM\Rightarrow 
CDAll_n(X)).
$$
This formula states that the matrix $M$ has the order~$n$  and commutes only with
scalar matrices from $D_n^Z(R)$. Therefore,
 $M=D_\rho S_\rho$, where $\rho$ is a cycle of length~$n$. Without loss
of generality, we can assume that $\rho=(1,2,\dots,n)$. It follows from 
Lemma~7 ([2]) that $M=\Phi_{N'}(S_\rho)$ for some $N'\in \Gamma_n(R)$.

Let us fix some matrix~$M$ satisfying the formula $Cycle_n(M)$ (if we do it, 
then a matrix $N'$ is chosen up to multiiplication by matrices
commuting with~$S_\rho$).

Let us consider now some matrix $M'$ satisfying the following formula
(with respect to~$M$):
\begin{multline*}
Trans_n(M,M'):=({M'}^2=1)\land ((MM')^{n-1}=1)
\land \exists X (CDOneMany_n(X)\land 
KOneMany_n(X,MM')).
\end{multline*}
Consider the matrix $\Phi_{N'}^{-1}(M')=D_\sigma S_\sigma$.

Since this natrix satisfies the condition $(D_\sigma S_\sigma)^2=1$, we have
$\sigma^2=1\Rightarrow \sigma=(i,j)$, $D_\sigma=diag[d_1,\dots,d_n]$,
$d_k=1$ for $k\ne i$, $k\ne j$, $d_id_j=1$.
Other conditions imply that the element $\sigma\rho$ is a cycle of
the length $n-1$. From the other hand 
$\sigma \rho=(i,j)(1,2,,\dots,n)=(1,2,\dots,i-1,j,j+1,\dots,n)
(i,i+1,\dots,j-1)$, therefore, $j=i+1$. Thus, $S_\sigma=S_{(i,i+1)}$,
$\sigma \rho=(1,2,\dots,i-1,i+1,\dots,n)$. From the condition
 $(MM')^{n-1}=1$ we have
$(diag[d_1,\dots,d_n] S_{(1,2,\dots,i-1,i+1,\dots,n)})^{n-1}=1$, therefore
$d_i^{n-1}=1 \Rightarrow d_i=1$. So $\Phi_{N'}^{-1}(M')=S_{(i,i+1)}$ 
for some $i=1,\dots,n$.

Consider then the matrix $N''=S_{\rho'}$, where $\rho'=\rho^{i-1}$. 
It is clear that
$\Phi_{N''}(S_\rho)=S_\rho$. We have here $\Phi_{N''}(S_{(i,i+1)})=S_{(1,2)}$. 
Therefore we have matrices $M_1$ and $M_2$, for which there exists 
 a matrix $N=N''N'$ such that $M_1=\Phi_N(S_{(1,2,\dots,n)})$, 
$M_2=\Phi_N (S_{(1,2)})$.

A matrix $M_1$ is an arbitrary matrix satisfying the formula $Cycle_n(M)$.

A matrix $M_2$  is an arbitrary matrix satisfying the formula $Trans_n(M_1,M)$.

The formula $Perm_n(M_1,M_2,M)$ can be constructed by the following:
for a given group $\Sigma_n$ 
for every its element~$\sigma$ we find  $\sigma=
(1,2)^{i_1}(1,2,\dots,n)^{j_1}\dots(1,2)^{i_k}(1,2,\dots,n)^{j_k}$.
Suppose that the elements of  $\Sigma_n$ are $\sigma_1,\dots,\sigma_N$,
$$
\sigma_l=(1,2)^{i_1^l}(1,\dots,n)^{j_1^l}\dots (1,2)^{i_{k(l)}^l}(1,\dots,n)^{j_{k(l)}^l}.
$$
Then
$$
Perm_n(M_1,M_2,M):= \bigvee_{l=1}^N (M=M_2^{i_1^l} M_1^{j_1^l}\dots M_2^{i_{k(l)}^l}
M_1^{j_{k(l)}^l}).
$$

For example, if $n=3$, we have
$\sigma_1=1$, $\sigma_2=(1,2,3)$, $\sigma_3=(1,2)$, $\sigma_4=(3,2,1)=(1,2,3)^2$,
$\sigma_5=(1,3)=(1,2,3)(1,2)$, $\sigma_6=(2,3)=(1,2)(1,2,3)$,
\begin{multline*}
Perm_3(M_1,M_2,M):= (M=1)\lor (M=M_1)\lor (M=M_2)\lor (M=M_1^2)
\lor (M=M_1M_2)\lor (M=M_2M_1).
\end{multline*}
\end{proof}

Now suppose that the matrices $M_1$ and $M_2$ are fixed. Therefore
the matrix~$N$ is fixed up to multiplication to a matrix $\alpha I\in D_n^Z(R)$.

\begin{lemma}\label{l9}
There exist formulas $GDOneMany_n(M_1,M_2,M)$ and
$DOneMany_n(M_1,M_2,M)$, which are true in the semigroup
 $G_n(R)$ if and only if $M$ has the form
 $\Phi_N(diag [\alpha,\beta,\dots,\beta])$, $\alpha,\beta\in
R_+^*$ and $\Phi_N(diag[\alpha,\beta,\dots,\beta])$, 
$\alpha,\beta\in R_+^*$, $\alpha\ne \beta$,
respectively.
\end{lemma}
\begin{proof}
It is clear that the obtained formulas are
$$
GDOneMany_n(M_1, M_2,M):=((M_2M_1)\cdot M=M\cdot (M_2 M_1))\land
Diag(M)
$$
and
$$
DOneMany_n(M_1,M_2,M):= GDOneMany_n(M_1,M_2,M)\land (M_2M\ne MM_2)
$$
(see Lemma~8 of the paper~[2]).
\end{proof}

\begin{lemma}\label{l10}
There exists a formula $G_2CD_{n-2}(M_1,M_2,M)$, which is true in the semigroup
 $G_n(R)$ if and only if
$M=\Phi_N(diag[X,a,\dots,a])$, where $X\in G_2(R)$, $a\in Z_+^*(R^*)$.
\end{lemma}
\begin{proof}
Similarly to the previous lemma we can write a formula
$DTwoMany_n(M_1,M_2,M)$, which is true in
$G_n(R)$ if and only if
$$
M=\Phi_N(diag[\alpha,\alpha,\beta,\dots,\beta])\in D_n(R),\quad \alpha\ne \beta.
$$

The formula
$$
DTransp_{1,2}(M_1,M_2,M):=(M^2=1)\land \exists (Diag(X)\land M=XM_2)
$$
defines involutions
$$
\Phi_N(diag[\xi,\xi^{-1},1,\dots,1]S_{(1,2)}),\quad \xi\in R_+^*
$$
(see the beginning of the proof of Lemma~9 of the paper~[2]).

Similarly , the formula
$$
CDTransp_{1,2}(M_1,M_2,M):=(M^2=1)\land \exists X (CDiad(X)\land M=XM_2)
$$
defines involutions
$$
\Phi_N(diag[\xi,\xi^{-1},1,\dots,1]S_{(1,2)}),\quad \xi\in Z_+^*(R^*).
$$

The formula
$$
CD_2D_{n-2}(M_1,M_2,M):=DTwoMany(M_1,M_2,M)
\land \forall X(DTransp_{
1,2}(M_1,M_2,X)\Rightarrow XM=MX)
$$
defines the set of matrices
$$
\Phi_N(diag[\mu,\mu,\eta,\dots,\eta]),\quad \mu\in Z_+^*(R^*),\eta \in R^*.
$$

As we remember, the formula $CDAll(M)$ defines matrices $\alpha I$, $\alpha
\in Z_+^*(R^*)$.

Now we will write the formula
\begin{multline*}
G_2CD_{n-2}(M_1,M_2,M):=\forall X (CD_2D_{n-2}(M_1,M_2,X)\Rightarrow\\
\Rightarrow
\exists Y(CDAll(Y)\land MXY=XYM\land\\
\land (M_1^{i_1}M_2^{j_2}\dots M_1^{i_k}M_2^{j_k})M=M(M_1^{i_1}M_2^{j_1}\dots M_1^{i_k}
M_2^{j_k}))),
\end{multline*}
where
$$
(1,2,\dots,n)^{i_1}(1,2)^{j_1}\dots (1,2,\dots,n)^{i_k}(1,2)^{j_k}=(3,\dots,n)
$$
(for $n=3$ we do not need the last condition).

This formula is equivalent to the assertion (2) from the proof
of Lemma~9 ([2]).
In the paper~[2] it was proved that in this case
$$
M=\Phi_N(diag[X,a,\dots,a]),\quad X\in G_2(R),\ a\in Z_+^*(R^*),
$$
what we needed to prove.
\end{proof}

\begin{lemma}\label{l11}
There exists a formula $ZDOneMany_n(M_1,M_2,M)$, which is true in
$G_n(R)$ if and only if
$$
M=\Phi_N(diag[\xi,\eta,\dots,\eta]),\quad \xi,\eta\in Z_+^*(R).
$$
\end{lemma}
\begin{proof}
Similarly to the formula $G_2CD_{n-2}(M_1,M_2,M)$ we can write the formula
$CD_{n-2}G_2(M_1,M_2,M)$, which is true if and only if
$$
M=\Phi_N(diag[a,\dots,a,X]),\quad X\in G_2(R),\ a\in Z_+^*(R).
$$

The formula 
$$
ZDAll(M):=\forall X (XM=MX)
$$
defines the center of $G_n(R)$, consisting of matrices $\alpha I$, $\alpha
\in Z_+^*(R)$.

Let for $n\ge l\ge 2$
$$
(1,l)=(1,2,\dots,n)^{i_1^l}(1,2)^{j_1^l}\dots (1,2,\dots,n)^{i_{k(l)}^l}
(1,2)^{j_{k(l)}^l}.
$$
Then the formula
\begin{multline*}
ZDOneMany_n(M_1,M_2,M):= DOneMany_n(M_1,M_2,M)
\land\\
\land \forall X (CD_{n-2}G_2(M_1,M_2,X)\Rightarrow MX=XM)\land \\
\land ZDAll_n(M\cdot (M_1^{i_1^2}M_2^{j_1^2}\dots M_1^{i_{k(2)}^2}
M_2^{j_{k(2)}^2}
\cdot MM_2^{-j_{k(2)}^2}M_1^{-i_{k(2)}^2}\dots M_2^{-j_1^2}M_1^{-i_1^2})\dots\\
\dots  (M_1^{i_1^n} M_2^{j_1^n}\dots M_1^{i_{k(n)}^n} M_2^{j_{k(n)}^n}
MM_2^{-j_{k(n)}^n} M_1^{-j_{k(n)}^n}\dots M_2^{-j_1^n}M_1^{-i_1^n}))
\end{multline*}
satisfies the condition (see the proof of Lemma~10 of~[2]).
\end{proof}

\begin{lemma}\label{l12}
There exists a formula $Main(M_1,M_2,M)$ which is true if and only if
either 
$$
M=\Phi_N B_{12}(x)=N B_{12}(x) N^{-1},
$$
or
$$
M=\Phi_N B_{21}(x)=NB_{21}(x) N^{-1}
$$
for some $x\in R_+$.
\end{lemma}

\begin{proof}
Consider the formula
\begin{multline*}
Main(M_1,M_2,M):= G_2CD_{n-2}(M_1,M_2,M)\land\\
\land \exists X(ZDOneMany
(M_1,M_2,X)
\land M^2=XMX^{-1})\land\\
\land \forall X (ZDOneMany(M_1,M_2,X)
\Rightarrow  M(XMX^{-1})=(XMX^{-1})M).
\end{multline*}
If a matrix $M$ satisfies the formula $G_2CD_{n-2}(M_1,M_2,M)$, then
$$
M=\Phi_N \begin{pmatrix}
\alpha& \beta& & &\\
\gamma& \delta& & &\\
& & a& &\\
& & &\ddots& \\
& & & & a
\end{pmatrix},\quad a\in Z_+^*(R^*),\ \begin{pmatrix}
\alpha& \beta\\
\gamma& \delta
\end{pmatrix}\in G_2(R).
$$
Since the matrix $M$ satisfies the subformula
$\forall X (ZDOneMany(M_1,M_2,X)\Rightarrow
M(XMX^{-1})=(XMX^{-1})M)$,
we have that
\begin{multline*}
\begin{pmatrix}
\alpha& \beta& & &\\
\gamma& \delta& & &\\
& & a& &\\
& & &\ddots& \\
& & & & a
\end{pmatrix}\begin{pmatrix}
\alpha& \xi \eta^{-1}\beta& & &\\
\eta \xi^{-1}\gamma& \delta& & &\\
& & a& &\\
& & &\ddots& \\
& & & & a
\end{pmatrix}=   \\ =
\begin{pmatrix}
\alpha& \xi \eta^{-1}\beta& & &\\
\eta\xi^{-1}\gamma& \delta& & &\\
& & a& &\\
& & &\ddots& \\
& & & & a
\end{pmatrix}\begin{pmatrix}
\alpha& \beta& & &\\
\gamma& \delta& & &\\
& & a& &\\
& & &\ddots& \\
& & & & a
\end{pmatrix}
\end{multline*}
for all $\eta,\xi\in Z_+^*(R)$, therefore
$$
\eta \xi^{-1}\beta\gamma=\xi \eta^{-1} \beta\gamma
$$
for all $\eta,\xi\in Z_+^*(R)$, and also for $\eta=2$, $\xi=1$.
So we have 
$$
2\beta\gamma=\frac{1}{2} \beta\gamma\Rightarrow \left( 2-\frac{1}{2}\right)
\beta\gamma=0.
$$
Consequently, either $\beta=0$, or $\gamma=0$.

The condition $\exists X (ZDOneMany(M_1,M_2,X)\land M^2=XMX^{-1})$,
$X=diag[\xi,\eta,\dots,\eta]$gives us either the relation
$$
\begin{pmatrix}
\alpha^2& \alpha\beta+\beta\delta& & &\\
0& \delta^2& & &\\
& & a^2& &\\
& & &\ddots& \\
& & & & a^2
\end{pmatrix}=
\begin{pmatrix}
\alpha& \xi\eta^{-1}\beta& & &\\
0& \delta& & &\\
& & a& &\\
& & &\ddots& \\
& & & & a
\end{pmatrix},
$$
and so $\alpha=\delta=a=1$, or the relation
$$
\begin{pmatrix}
\alpha^2& 0& & &\\
\alpha \gamma+\gamma \delta& \delta^2& & &\\
& & a^2& &\\
& & &\ddots& \\
& & & & a^2
\end{pmatrix}
\begin{pmatrix}
\alpha& 0& & &\\
\eta \xi^{-1}\gamma& \delta& & &\\
& & a& &\\
& & &\ddots& \\
& & & & a
\end{pmatrix},
$$
and so $\alpha=\delta=a=1$.
Therefore, either $M=B_{12}(\beta)$, or $M=B_{21}(\gamma)$.
\end{proof}

\begin{lemma}\label{l13}
There exists a formula
$$
MainUnit_{1,2}(M_1,M_2,M),
$$
which is true only for the matrix $M=\Phi_N B_{12}(1)$.
\end{lemma}
\begin{proof}
Consider the formula
\begin{multline*}
MainUnit_{1,2}(M_1,M_2,M)=Main(M_1,M_2,M)\land (M\ne 1)\land\\
\land \bigl((M_1^{i_1^1}M_2^{j_1^1}\dots M_1^{i_{k(1)}^1}M_2^{j_{k(1)}^1})
\cdot M\cdot 
(M_2^{-j_{k(1)}^1}M_1^{-i_{k(1)}^1}\dots M_2^{-j_1^1}M_1^{-i_1^1})\cdot\\
\cdot (M_1^{i_1^2}M_2^{j_1^2}\dots M_1^{i_{k(2)}^2} M_2^{j_{k(2)}^2})
\cdot M\cdot 
(M_1^{i_1^2} M_2^{j_1^2}\dots M_1^{i_{k(2)}^2}M_2^{j_{k(2)}^2})=\\
=(M_1^{i_1^2}M_2^{j_1^2}\dots M_1^{i_{k(2)}^2}M_2^{j_{k(2)}^2})\cdot M\cdot\\
\cdot (M_1^{i_1^2} M_2^{j_1^2}\dots M_1^{i_{k(2)}^2}M_2^{j_{k(2)}^2})
(M_1^{i_1^1}
M_2^{j_1^1}\dots M_1^{i_{k(1)}^1} M_2^{j_{k(1)}^1})\cdot M\cdot\\
\cdot (M_1^{i_1^1} M_2^{j_1^1} \dots M_1^{i_{k(1)}^1} 
M_2^{j_{k(1)}^1})\cdot M\bigr),
\end{multline*}
where
\begin{align*}
(2,3)& = (1,2,\dots,n)^{i_1^1}(1,2)^{j_1^1}\dots (1,2,\dots,n)^{i_{k(1)}^1} 
(1,2)^{j_{k(1)}^1},\\
(1,3)& = (1,2,\dots,n)^{i_1^2}(1,2)^{j_1^2}\dots (1,2,\dots,n)^{i_{k(2)}^2} (1,2)^{
j_{k(2)}^2}.
\end{align*}

Since every matrix satisfying the formula $Main(M_1,M_2, M)$, can have 
either the form $\Phi_N B_{12}(x)$, or the form $\Phi_N B_{21}(x)$, 
we have one of two conditions: either
\begin{multline*}
S_{(2,3)} \Phi_N^{-1}(M) S_{(2,3)} S_{(1,3)} \Phi_N^{-1}(M)
 S_{(1,3)}= S_{(1,3)} \Phi_N^{-1}(M) S_{(1,3)} S_{(2,3)}
\Phi_N^{-1}(M) S_{(2,3)} \Phi_N^{-1}(M)\Rightarrow\\
\Rightarrow
\begin{pmatrix}
1& x^2& x\\
0& 1& 0\\
0& x& 1
\end{pmatrix} =
\begin{pmatrix}
1& x& x\\
0& 1& 0\\
0& x& 1
\end{pmatrix}\Rightarrow\\
x^2=x\Rightarrow x=1\Rightarrow M=\Phi_N B_{12}(1),
\end{multline*}
or
\begin{multline*}
S_{(2,3)} \Phi_N^{-1}(M) S_{(2,3)} S_{(1,3)} \Phi_N^{-1}(M)
 S_{(1,3)}= S_{(1,3)} \Phi_N^{-1}(M) S_{(1,3)} S_{(2,3)}
\Phi_N^{-1}(M) S_{(2,3)} \Phi_N^{-1}(M)\Rightarrow\\
\Rightarrow
\begin{pmatrix}
1& 0& 0\\
0& 1& x\\
x& 0& 1
\end{pmatrix} =
\begin{pmatrix}
1& 0& 0\\
x^2+x& 1& x\\
x& 0& 1
\end{pmatrix}\Rightarrow
x^2+x=0,
\end{multline*}
but it is impossible.
\end{proof}

\begin{lemma}\label{l14}
There exists a formula $Main_{1,2}(M_1,M_2,M)$, which is true in $G_n(R)$
for the matrices $\Phi_N (B_{12}(x))$,
$x\in R\cup \{ 0\}$, and only for them.
\end{lemma}

\begin{proof}
Consider the formula
$$
Main_{1,2}(M_1,M_2,M):= Main(M_1,M_2,M)
\land \forall X(MainUnit_{1,2}(X)\Rightarrow XM=MX).
$$
Since a matrix $M$ satisfies the formula $Main(M_1,M_2,M)$, we have that
either $M=\Phi_N(B_{12}(x))$ or $\Phi_N(B_{21}(x))$ for some $x\in R$.
In the second case the second part of the formula implies
$$
\begin{pmatrix}
1& 1\\
0& 1
\end{pmatrix} \begin{pmatrix}
1& 0\\
x& 1
\end{pmatrix}=\begin{pmatrix}
1& 0\\
x& 1
\end{pmatrix}
\begin{pmatrix}
1& 1\\
0& 1
\end{pmatrix}
\Rightarrow 
\begin{pmatrix}
1+x& 1\\
x& 1
\end{pmatrix}=
\begin{pmatrix} 
1& 1\\
x& 1+x
\end{pmatrix}\Rightarrow x=0.
$$
\end{proof}
\begin{lemma}\label{l15}
There exist formulas
$$
Addit_n(M_1,M_2,X_1,X_2,X_3)\text{ and }
Multipl_n(M_1,M_2,X_1,X_2,X_3),
$$
which are true in $G_n(R)$ if and only if
 $X_1=\Phi_N (B_{12}(x_1))$, $X_2=\Phi_N(
B_{12}(x_2))$, $X_3=\Phi_N (B_{12}(x_3))$, 
where $x_1,x_2,x_3\in R_+\cup \{ 0\}$, 
and, respectively, either $x_3=x_1+x_2$, or $x_3=x_1\cdot x_2$.
\end{lemma}

\begin{proof}
The obtained formulas are
\begin{multline*}
Addit_n(M_1,M_2,X_1,X_2,X_3):= Main_{1,2}(M_1,M_2,X_1)\land\\
\land Main_{1,2}(M_1,M_2,X_2)\land Main_{1,2} (M_1,M_2,X_3)\land
X_3=X_1\cdot X_2
\end{multline*}
and
\begin{multline*}
Multipl_n(M_1,M_2,X_1,X_2,X_3):= Main_{1,2}(M_1,M_2,X_1)\land\\
\land Main_{1,2}(M_1,M_2,X_2)\land Main_{1,2}(M_1,M_2,X_3)\land\\
\land (M_1^{i_1^1}M_2^{j_1^1}\dots M_1^{i_{k(1)}^1} M_2^{j_{k(1)}^1})\cdot X_1\
\cdot (M_1^{i_1^1}M_2^{j_1^1}\dots M_1^{i_{k(1)}^1} M_2^{j_{k(1)}^1})\times\\ 
\times (M_1^{i_1^2}M_2^{j_1^2}\dots M_1^{i_{k(2)}^2} M_2^{j_{k(2)}^2})
\cdot X_2\cdot (M_1^{i_1^2}M_2^{j_1^2}\dots M_1^{i_{k(2)}^2} M_2^{j_{k(2)}^2})
=\\
=(M_1^{i_1^2}M_2^{j_1^2}\dots M_1^{i_{k(2)}^2} M_2^{j_{k(2)}^2})\cdot X_2 \cdot 
(M_1^{i_1^2}M_2^{j_1^2}\dots M_1^{i_{k(2)}^2} M_2^{j_{k(2)}^2})\times\\
\times (M_1^{i_1^1}M_2^{j_1^1}\dots M_1^{i_{k(1)}^1} M_2^{j_{k(1)}^1})\cdot X_1
\cdot (M_1^{i_1^1}M_2^{j_1^1}\dots M_1^{i_{k(1)}^1} M_2^{j_{k(1)}^1})\cdot X_3,
\end{multline*}
where 
\begin{align*}
(2,3)& = (1,2,\dots,n)^{i_1^1}(1,2)^{j_1^1}\dots (1,2,\dots,n)^{i_{k(1)}^1} 
(1,2)^{j_{k(1)}^1},\\
(1,3)& = (1,2,\dots,n)^{i_1^2}(1,2)^{j_1^2}\dots (1,2,\dots,n)^{i_{k(2)}^2} (1,2)^{
j_{k(2)}^2}.
\end{align*}

The first formula is absolutely clear, and the second one 
follows from the relation
$$
B_{13}(x_1)B_{32}(x_2)=B_{32}(x_2)B_{13}(x_1)B_{12}(x_1x_2).
$$
\end{proof}

\setcounter{theorem}{0}
\begin{theorem}
Let $R$ and $S$ be linearly ordered rings with~$1/2$, $G_n(R)$ and $G_m(S)$ 
be semigroups of invertible matrices with nonnegative elements, $n,m\ge 3$. 
Then the semigroups $G_n(R)$ and $G_m(S)$ are elementary equivalent if and
only if $n=m$
and the semirigns $R_+$ and $S_+$ are elementary equivalent.
\end{theorem}

\begin{proof}
1. Let   $R_+\equiv S_+$.
Let us show that $G_n(R)\equiv G_n(S)$.

Actually, let us consider the semigroup $G_n(R)$
($n\ge 3$) and some sentence of group language
$$
{\cal U}=(Q_1 X_1)\dots (Q_r X_r){\cal B}(X_1,\dots,X_r)\
(Q_i=\exists, \forall).
$$

A matrix $X\in G_n(R)$ is a set of $n^2$ elements of the semiring
$R_+\cup \{ 0\}$
$\{ x_{ij}\}$ with the condition that in $GL_n(R)$ there exists 
an inverse matrix. This condition can be expressed by the following. 
Let ${\cal J}_n$ be the set $\{ 1,\dots,n\}$. Then we have
\begin{multline*}
\exists y_{11},\dots,\exists y_{nn}:
\bigvee_{S\subseteq {\cal J}_n\times {\cal J}_n} 
\bigwedge_{i,k\in {\cal J}_n\times {\cal J}_n}
\Bigl(\delta_{ik}+\sum_{\langle j,k\rangle \in S}x_{ij}y_{jk}=
\sum_{\langle j,k\rangle\notin S} 
x_{ij}y_{jk} \bigwedge\\
\bigwedge \delta_{ik}+\sum_{\langle i,j\rangle
\in S}y_{ij}x_{ij}=\sum_{\langle i,j\rangle\notin 
S} y_{ij}x_{jk}\Bigr).
\end{multline*}

Let us denote this condition (for a matrix $X$) by ${\mathbf G}(X)$. 

Then  $\cal U$ is true in $G_n(R)$ if and only if the following sentence
 ${\cal U}_R$ is true in~$R_+$: it is obtained from~$\cal U$ 
by the following process
(compare with~[3]):

 ) in the formula ${\cal B}(X_1,\dots ,X_r)$ all relations $X_i=X_j$ and
$X_i=X_jX_k$ are changed respectvely to the formulas
$$
\bigwedge_{\lambda,\mu}(x_{i\lambda\mu}=x_{j\lambda\mu})\text{ and }
\bigwedge_{\lambda,\mu}(x_{i\lambda\mu}=x_{j\lambda1}x_{k1\mu}+\dots+
x_{j\lambda n}x_{kn\mu}).
$$

¡) If the formula ${\cal B}_{i+1}=(Q_{i+1}X_{i+1})\dots (Q_r X_r){\cal B}$
is translated to the formula ${\cal Q}_{i+1}$, then the formula
 ${\cal B}_i=
(\forall X_i){\cal B}_{i+1}$ is translated to the formula
$$
\forall x_{i11}\dots \forall x_{inn}(\mathbf G(x_{i11},\dots,x_{inn})\Rightarrow 
{\cal Q}_{i+1}),
$$
and the formula ${\cal B}_i=(\exists X_i){\cal B}_{i+1}$ is translated 
to the formula
$$
\exists x_{i11}\dots \exists x_{inn}
(\mathbf G(x_{i11},\dots, x_{inn})\land {\cal Q}_{i+1}).
$$

Consequently, every sentence~$\cal U$ of group language
can be (by this method) translated
to the sentence~${\cal U}_R$ of ring language.
The form of~${\cal U}_R$ does not depend of a basic field.
So  $\cal U$ is true in $G_n(R)$ if and only if
${\cal U}_R$ is true in~$R$. We know that ${\cal U}_R$ is true in~$R$ 
if and only if ${\cal U}_R$ is true~$S$, since $R\equiv S$.
Therefore $\cal U$ is true in $G_n(R)$
if and only if ${\cal U}$ is true in $G_n(S)$, and, consequently,
$G_n(R)\equiv G_n(S)$.

2. Now we will prove the converse implication.
Let  semigroups $G_n(R)$ and
$G_m(S)$ be elementary equivalent. Lemma~\ref{l4} implies 
$m=n$, therefore we can suppose that we have semigroups $G_n(R)$ and
$G_n(S)$.

Now we will prove that $R_+\equiv S_+$. Suppose that we have some
sentence $\varphi$ of ring language. Let us translate it to the sentence
 $\overline \varphi$
of group language by the following algorhitm:

the sentence $\overline \varphi$ has the form
$$
\exists M_1 \exists M_2 (Cycle_n(M_1)\land Transp_n(M_1,M_2)\land
\varphi' (M_1,M_2)),
$$
where the formula $\varphi$ is obtained from the sentence $\varphi$ 
by the following translations of subformulas of~$\varphi$:

1) the subformula $\forall x \psi$ is translated to the subformula
 $\forall X (Main_{1,2}
(M_1,M_2,X)\Rightarrow \psi')$;

2) the subformula $\exists x \psi$ is translated to the subformula
 $\exists X (Main_{1,2}
(M_1,M_2,X)\land \psi')$;

3) the subformula $x = y$ is translated to the subformula $X = Y$;

4) the subformula $x = y+z$ is translated to the subformula
 $Addit_n(M_1,M_2,X,Y,Z)$;

5) the subformula $x = yz$ is translated to the subformula
 $Multipl_n(M_1,M_2,X,Y,Z)$.

It is clear that the sentence $\varphi$ is true in the semiring~$R_+$ 
if and only if the sentence $\overline \varphi$ is true in the semigroup
 $G_n(R)$.
Since the semigroups $G_n(R)$ and $G_n(S)$ are elementary equivalent, we have that
the sentence $\overline \varphi$ is true in $G_n(R)$
if and only if it is true in $G_n(S)$, and it is equivalent to 
 $\varphi$ being true in~$S_+$.
Therefore the semirings $R_+$ and $S_+$ are elementary equivalent.

\end{proof}


\begin{thebibliography}{99}

\bibitem{1} A.\,V. Mikhalev, M.\,A. Shatalova. Automorphisms and
antiautomorphisms of the semigroup of invertible matrices with
nonnegative elements.  Mat. Sbornik, 81(4), 1970, 600--609.

\bibitem{2} E.I. Bunina, A.V. Mikhalev. Automorphisms of the semigroup
of invertible matrices with nonnegative elements. 
Fund. i Prikl. Mat., 11(2), 2005, 3--23.


\bibitem{3} A.I. Maltsev. On elementary properties of linear groups. 
Problems of Mathematics and Mechanics, Novosibirsk, 1961, 110--132.

\bibitem{4} C.I. Beidar, A.V. Mikhalev. On Malcev's theorem on elementary 
equivalence of linear groups. Contemporary mathematics, 131, 1992, 
29--35

\bibitem{5} E.I. Bunina. Elementary equivalence of unitary linear groups over fields.
Fund. i Prikl. Mat., 4(4), 1998, 1265--1278.

\bibitem{6} E.I. Bunina. Elementary equivalence of unitary linear groups over rings
and skewfields. Russian Math. Surveys, 53(2), 1998, 137--138.

\bibitem{7} E.I. Bunina. Chevalley groups over fields and their elementary properties.
Russian Math. Surveys, 59(5), 2004, 952--953.
\end{thebibliography}
\end{document}